\documentclass[a4paper]{amsart}
\usepackage{amssymb}
\usepackage{latexsym}
\usepackage{euscript}
\usepackage{hyperref}
\def\sgn{\mathop{\rm sgn}\nolimits}
\def\tr{\mathop{\rm tr}\nolimits}

\newtheorem{theorem}{Theorem}
\theoremstyle{plain}

\newtheorem{definition}{Definition}

\newtheorem{lemma}{Lemma}

\newtheorem{proposition}{Proposition}
\newtheorem{remark}{Remark}

\numberwithin{equation}{section}

\def\grad{\mathop{\rm \pmb{grad}}\nolimits}
\def\Hess{\mathop{\rm Hess}\nolimits}

\newcommand{\xig}{\pmb{\xi}}
\newcommand{\zg}{\pmb{\zeta}}
\newcommand{\eg}{\pmb{\eta}}
\newcommand{\nog}{\pmb{n}}
\newcommand{\ug}{\pmb{u}}
\newcommand{\vg}{\pmb{v}}
\newcommand{\sg}{\pmb{\sigma}}

\newcommand{\nablag}{\pmb{\nabla}}

\newcommand{\R}{\mathbb{R}}

\newcommand{\ES}{\EuScript{S}}
\newcommand{\EL}{\EuScript{L}}
\newcommand{\EH}{\EuScript{H}}
\newcommand{\EHb}{\overline{\EuScript{H}}}

\newcommand{\dnu}[1]{%
\dfrac{\partial #1}{\partial n}
}

\newcommand{\DOT}[2]{
\left\langle #1,#2 \right\rangle}

\begin{document}
\title{Generalized Kirchhoff approximation for Helmholtz equation}
\author{F. Cuvelier}
\address{F.Cuvelier : Universit\'e Paris XIII, Institut Galil\'ee, LAGA CNRS UMR
7539, 93430 Villetaneuse, {\textsc{France}} }
\email{cuvelier@math.univ-paris13.fr}

\begin{abstract}
We give integral formulas to approximate solutions of Dirichlet and Neumann
problems for Helmholtz equation at high frequencies. These approximations
are valid in the complementary of a union of convex compact obstacles. The
first step of the iterative procedure is the classical Kirchhoff
approximation. Convergence is proved by comparison with the geometrical
optics asymptotics. The method is shown to be numerically stable.
\end{abstract}

\maketitle

\section{Introduction}

Let $\Omega $ be an open set in $\mathbb{R}^{3}$ and $\Omega ^{\prime }=%
\mathbb{R}^{3}\setminus \Omega $. We study the high frequency diffraction
problems of an incident plane wave in $\Omega $ for Helmholtz equation with
respectively, Dirichlet and Neumann boundary conditions :
\begin{equation}
\left\{
\begin{array}{rcl}
\Delta v(x)+k^{2}v(x) & = & 0\ \ \ \mbox{for}\ x\in \Omega , \\
v(x) & = & 0\ \ \ \mbox{for}\ x\in \partial \Omega , \\
v(x) & = & e^{-ik\DOT{\xi}{x} }+u(x)\ \ \ \mbox{for}\
x\in \Omega ,
\end{array}
\right.  \tag{$D$}  \label{Dirichlet}
\end{equation}
\begin{equation}
\left\{
\begin{array}{rcl}
\Delta v(x)+k^{2}v(x) & = & 0\ \ \ \mbox{for}\ x\in \Omega , \\
\dnu{v}(x) & = & 0\ \ \ \mbox{for}\ x\in \partial
\Omega , \\
v(x) & = & e^{-ik\DOT{\xi}{x}}+u(x)\ \ \ \mbox{for}\
x\in \Omega .
\end{array}
\right.  \tag{$N$}  \label{Neumann}
\end{equation}
Here $u$ satisfies the Sommerfeld radiation condition :
\begin{equation*}
r^{2}({\frac{\partial u}{\partial r}}+iku)\mbox{ bounded\ when }r=\mid x\mid
\rightarrow +\infty .
\end{equation*}
The incident plane wave, $e^{-ik\DOT{\xi}{x}},$ is given with
 the normalization $\mid \xi \mid =1$. Here, \textit{high frequency}
means that the wave length is small with respect to $\partial\Omega$ curvatures. 
So usual numerical methods, such as finite element method, boundary
element method and so on, fall down. 

A classical high frequency approximation is given by geometrical optics ,
which, for a point $x\in \Omega $, allows us to compute, from optic rays
going through $x$, an approximation of the diffracted wave by $\Omega
^{\prime }$. We obtain for Dirichlet problem $(D)$ and Neumann problem $(N)$
respectively
\begin{equation}
v_{O.G.}^{D}(x)=\sum e^{-ik\varphi (x)}a_{0}^{D}(x)
\label{Geometrical Optic D}
\end{equation}
and
\begin{equation}
v_{O.G.}^{N}(x)=\sum e^{-ik\varphi (x)}a_{0}^{N}(x).
\label{Geometrical Optic N}
\end{equation}
Here, the phase $\varphi (x)$ is the length of optic rays going through $x$
and computation of the amplitude $a_{0}^{\cdot }(x)$ work out by propagation
and reflection formulas along optic rays (see \cite{Cuv:GO:2013}). The
main problem of this method is its numerical instability : in order to compute 
this approximation, it's necessary to determine \textbf{all} the optic rays 
going through $x.$ But, small errors in the numerical representation of $\partial
\Omega $ can give large errors in the optic rays determination.

An other one is Kirchhoff approximation, based on integral
representations. We give Kirchhoff approximation respectively for Dirichlet
problem $(D)$ and Neumann problem $(N)$ :
\begin{equation}
v^D_{\mbox{Kir.}}(x)=e^{-ik\DOT{\xig}{x}}+\frac{1}{4\pi }\int_{\partial
\Omega }ik\left( \DOT{\xig }{\nog(\sigma )} 
                - \left| \DOT{\xig}{\nog(\sigma )} \right| 
          \right) e^{-ik\DOT{\xig}{\sg}}
{\frac{e^{-ik\mid x-\sigma \mid }}{\mid x-\sigma \mid }}d\sigma  \label{Kirchhoff for Dirichlet}
\end{equation}
\begin{equation}
\begin{array}{c}
v^N_{\mbox{Kir.}}(x)
\\=\\
e^{-ik\DOT{\xig }{x} } \\
+\\ \frac{1}{4\pi }\int_{\partial \Omega }
ik\left( \frac{\DOT{\xig }{\nog(\sigma )} }
              {\left| \DOT{\xig }{\nog(\sigma )} \right| }
         -1\right) 
\DOT{\frac{x-\sigma }{\left| x-\sigma \right| }}{\nog(\sigma )}
e^{-ik\DOT{\xig}{\sg} }{
\frac{e^{-ik\mid x-\sigma \mid }}{\left| x-\sigma \right| }}d\sigma
\end{array}
\end{equation}
Here, $\nog(\sigma )$ is the unit normal to $\Gamma$ at point $
\sigma $, exterior to $\Omega ^{\prime }$. But, the validity of this method
is restricted to $\Omega ^{\prime }$ beeing a strictly convex compact (see \cite{Melrose:Taylor} 
for Dirichlet problem and \cite{Yingst79}, \cite{Yingst83}
for Neumann problem) and false otherwise. This is due to the incapacity of
this method to \textit{{see} }multiple reflections.

The purpose of this paper is to determine an iterative integral method,
numerically stable, equivalent, at first order and high frequency, to the
geometrical optic approximation sets $\Omega ^{\prime }$ is a finite and
disjointed union of regular and strictly convex compacts. In both case, the
first step is given by Kirchhoff approximation.

It relies in ...(écrire les principales étapes du papier)

\section{Notations and definitions}

\subsection{Gradient and Hessian on Surfaces}
Let $K\subset \mathbb{R}^3$ be a compact and $\Gamma$ it boundary. Let us suppose that $\Gamma$
is a regular and orientable surface.

\begin{definition}
(\textbf{Gradient on Surfaces}). The \textbf{gradient} of a differentiable function 
$\varphi : \Gamma\subset \mathbb{R}^3 \rightarrow \mathbb{R}$ is
a differentiable map $\grad \varphi : \Gamma \rightarrow \mathbb{R}^3$ which assigns to each point 
$\sigma\in \Gamma$ a
vector $\grad \varphi(\sigma)\in T_\sigma(\Gamma)\subset \mathbb{R}^3$ such that
$$<\grad \varphi(\sigma),v>_\sigma=d\varphi_\sigma(v)\ \ \ \forall v\in T_\sigma(\Gamma).$$
\end{definition}
\begin{definition}(\textbf{Hessian on Surfaces}).
The \textbf{hessian} of a twice differentiable function 
$\varphi : \Gamma\subset \mathbb{R}^3 \rightarrow \mathbb{R}$ is the function
$\Hess\varphi : \Gamma \rightarrow \EL(\mathbb{R^3})$ which assigns to 
each point $\sigma\in \Gamma$ a matrix 
$\Hess\varphi(\sigma)\in\EL\left(T_{\sigma}(\Gamma) \right)$ such that 
$$\DOT{\Hess \varphi(\sigma) v}{v}=d^2\varphi_\sigma(v).v\ \ \ \forall v\in T_\sigma(\Gamma).$$
\end{definition}
\begin{proposition}
With previous definition, and by taylor's expansion we obtain for 
$t\in \R$ and $v\in T_\sigma(\Gamma)$
$$\varphi(\sigma+t v) - \varphi(\sigma) = 
t \DOT{\grad \varphi(\sigma)}{v}
+\frac{t^2}{2}\DOT{\Hess \varphi(\sigma) v}{v} +o(t^2)$$
\end{proposition}

\begin{definition}
The gradient of a differentiable function
$\varphi : \overbrace{\Gamma\times \cdots \times \Gamma}^{N \mbox{ times}} \subset \R^{3N}
\rightarrow \R$ is a differentiable map 
$\grad \varphi : \Gamma\times \cdots \times \Gamma \rightarrow \R^{3N}$ which
assigns to each point 
$\nu=(\sigma_1,\ldots,\sigma_N)\in\Gamma^N$ a vector 
$\grad \varphi(\nu)\in \left(T_{\sigma_1}(\Gamma)\times\cdots\times T_{\sigma_N}(\Gamma)\right)
\subset \R^{3N}$ such that
$$\DOT{\grad \varphi(\nu)}{\varpi}_\nu=d\varphi_\nu(\varpi)\ \ \forall \varpi\in
\left(T_{\sigma_1}(\Gamma)\times\cdots\times T_{\sigma_N}(\Gamma)\right)$$
The function $\nablag_{\sigma_i} \varphi : \Gamma^N \rightarrow \R^3$
which assigns to each point 
$\nu=(\sigma_1,\ldots,\sigma_N)\in\Gamma^N$ a vector $\nablag_{\sigma_i} \varphi(\nu)\in
T_{\sigma_i}(\Gamma)\subset \R^3$ is defined by
$$\DOT{\nablag_{\sigma_i} \varphi(\nu)}{\varpi_i}=d\varphi_\nu(\varpi)\ \ \forall \varpi_i
\in  T_{\sigma_i}(\Gamma)$$
with $\varpi=(0,\ldots,0,\varpi_i,0,\ldots,0).$
\end{definition}
\begin{definition}
The \textbf{hessian} of a twice differentiable function $\varphi : \Gamma^N \rightarrow \R$
is the function 
$\Hess \varphi : \overbrace{\Gamma\times \cdots \times \Gamma}^{N \mbox{ times}} \rightarrow \EL\left(\R^{3N}\right)$ which assigns to
each point $\nu=(\sigma_1,\ldots,\sigma_N)\in\Gamma^N,$ a matrix $\Hess \varphi (\nu)\in 
\EL\left(T_{\sigma_1}(\Gamma)\times\cdots\times T_{\sigma_N}(\Gamma)\right)$ such that
$$\DOT{\Hess \varphi (\nu)\varpi}{\varpi}_\nu= d^2\varphi_\nu(\varpi).\varpi\ \ \forall
\varpi\in T_{\sigma_1}(\Gamma)\times\cdots\times T_{\sigma_N}(\Gamma).$$

The function $\EH_{i,j} \varphi : \overbrace{\Gamma\times \cdots \times \Gamma}^{N\mbox{ times}} \rightarrow
\EL\left(\R^3\right),$ $1\leq i\neq j \leq N,$ which assigns to each point 
$\nu=(\sigma_1,\ldots,\sigma_N)\in\Gamma^N,$ a matrix 
$\EH_{i,j} \varphi(\nu)\in \EL\left(T_{\sigma_i}(\Gamma),T_{\sigma_j}(\Gamma) \right)$
is defined by
$$\DOT{\EH_{i,j} \varphi(\nu)\omega_i}{\omega_j}=
d^2\varphi_\nu(\varpi).\varpi\ \ \forall (\omega_i,\omega_j)\in T_{\sigma_i}(\Gamma)\times
T_{\sigma_j}(\Gamma)$$
with $\varpi=(\varpi_1,\ldots,\varpi_N),$ $\varpi_k=0$ for $k\neq i$ and $k\neq j,$
$\varpi_i=\omega_i$ and $\varpi_j=\omega_j.$

The function $\EH_{i,i} \varphi : \overbrace{\Gamma\times \cdots \times \Gamma}^{N\mbox{ times}} \rightarrow
\EL\left(\R^3\right),$ $1\leq i \leq N,$ which assigns to each point 
$\nu=(\sigma_1,\ldots,\sigma_N)\in\Gamma^N,$ a matrix 
$\EH_{i,i} \varphi(\nu)\in \EL\left(T_{\sigma_i}(\Gamma)\right)$
is defined by
$$\DOT{\EH_{i,i} \varphi(\nu)\omega_i}{\omega_i}=
d^2\varphi_\nu(\varpi).\varpi\ \ \forall \omega_i\in T_{\sigma_i}(\Gamma)$$
with $\varpi=(\varpi_1,\ldots,\varpi_N)\in 
T_{\sigma_1}(\Gamma)\times\cdots\times T_{\sigma_N}(\Gamma),$ $\varpi_k=0$ for $k\neq i,$ and
$\varpi_i=\omega_i.$

\end{definition}

\subsection{Geometrical notations}
\begin{itemize}
\item Let $(K_{i})_{i=1,\cdots ,N}$ be a set of regular, disjoint and strictly convex
compact in $\mathbb{R}{^{3}.}$
\item We denote by $\Gamma_i,$ the boundary of $K_i,$ and 
$\Gamma =\bigcup_{i=1}^N \Gamma_i.$ \\
\item Thus $\Gamma_i$ is a regular and orientable surface. 
So, given a point $\sigma$ of surface $\Gamma_i$ we can 
choose the coordinate axis of $\mathbb{R}^3$ so that origin $O$ of the coordinates is at $\sigma$ and
the $z$ axis is directed along the \textbf{negative normal} (i.e. the outer normal) 
$\nog(\sigma)$ of $\Gamma_i$ in $\sigma$ (thus, the $xy$ plane
agrees with $T_\sigma(\Gamma_i)$ : tangent plane of $\Gamma_i$ in $\sigma$). It follows that a
neighborhooh of $\sigma$ in $\Gamma_i$ can be represented in the form 
$z=g_i(u,v),\ (u,v)\in U\subset \mathbb{R}^2,$ where $U$ is an open set and $g_i$ is a differentiable 
function with $g_i(0,0)=\frac{\partial g_i}{\partial u}(0,0)=\frac{\partial g_i}{\partial v}(0,0)=0.$

Let us assume further that the $u$ and $v$ axes are directed along the 
\textbf{principal directions}, 
with the axis $u$ along the direction of maximum principal curvature. Thus 
$$k_1^i(\sigma)=\frac{\partial^2 g_i}{\partial u^2}(0,0),\ 
k_2^i(\sigma)=\frac{\partial^2 g_i}{\partial v^2}(0,0),\ 
\mbox{ and } \frac{\partial^2 g_i}{\partial u\partial v}(0,0)=0$$
and, so we obtain by developing $g_i(u,v)$ into Taylor's expansion about $(0,0)$
$$g_i(u,v)=-\frac{1}{2}(k_1^i u^2 + k_2^i v^2) +o(u^2+v^2)$$
\item We note 
$R^i_u(\sigma)=\frac{1}{k_1^i(\sigma)}$ and
$R^i_v(\sigma)=\frac{1}{k_2^i(\sigma)}$ the \textbf{principal radii of curvature}.

\item Let us denote $\Re_i(\sigma)$ the orthonormal basis 
$\Re_i(\sigma)=\{ \ug_i, \vg_i, \nog_i\}$ where 
$\nog_i$ is the negative normal of $\Gamma_i$ in $\sigma,$
$\ug_i$ and $\vg_i$ are the principal directions of $\Gamma_i$ in $\sigma$ 
with $\ug_i$ the direction of maximum principal curvature. 

\item  We set $\Omega ^{\prime }=\bigcup\limits_{i=1}^{N}K_{i}$
and we note $\Omega _{i}=\mathbb{R}^3\setminus K_{i}.$





\item  Let $\mathcal{M}_{m,n}(\mathbb{R})$ the set of real matrix of size $%
m\times n.$

\item  We note
\begin{equation*}
\Gamma_{*}^{l}=\left\{ (\sigma _{1},\cdots ,\sigma _{l})\in
\Gamma^l \mbox{ with }\left( \sigma _{j}\in \Gamma_{p}\mbox{
and }\sigma _{j+1}\in \Gamma_{q}\Rightarrow p\neq q\right) \right\}
\end{equation*}
and the phase function $\psi _{l} : \R^3 \times \Gamma_*^l \rightarrow \R$ 
defined, for all $\nu=\sigma _{1},\cdots ,\sigma _{l}\in\Gamma_*^l,$ by
\begin{equation*}
\psi _{l}(x;\nu)=\DOT{\xig}{\sigma_1}
+\sum_{j=1}^{l-1}{\left| \sigma _{j+1}-\sigma _{j}\right| }%
+\left| x-\sigma _{l}\right|
\end{equation*}

\item  We call $\mathcal{C}_{l}(x)$ the set of all $l$-uplet $(\sigma _{1},\cdots ,\sigma _{l})
\in \Gamma_*^l$ such that:

\begin{enumerate}
\item  $\DOT{\xig}{\nog(\sigma_1)}<0$ and 
$\DOT{\sigma _{j+1}-\sigma _{j}}{\nog(\sigma _{j+1})}<0$ for all $j\in
\{1,\cdots ,l-1\}$,

\item  The phase $\psi _{l}(x,\bullet )$ is stationary on $\Gamma_*^l$ at point 
$\nu=(\sigma _{1},\cdots ,\sigma _{l})$ (i.e. $\grad_\nu \psi_l(x;\nu)=0$)
\end{enumerate}

We note $\mathcal{C}(x)=\bigcup_{l}\mathcal{C}_{l}(x)$.

\item  Let $\nu =(\sigma _{1}^{\nu },\cdots ,\sigma _{l}^{\nu })\in \mathcal{%
C}_{l}(x)$. 
We note $\Re _{j}^{\nu }=\Re_j(\sigma_j^\nu)=\{\ug_j^\nu,\vg_j^\nu,\nog_j^\nu\},\ \
j\in \{1,\cdots ,l\}.$

\item  Let $\nu =(\sigma _{1}^{\nu },\cdots ,\sigma _{l}^{\nu })\in \mathcal{%
C}_{l}(x),$ we define $\xig_{j}^{\nu }$, for $j\in \{1,\cdots ,l\}$ by :
\begin{equation*}
\sigma _{j+1}^{\nu }-\sigma _{j}^{\nu }=\lambda _{j}^{\nu }\xig_{j}^{\nu }\
\mbox{with }\ \lambda _{j}^{\nu }=\mid \sigma _{j+1}^{\nu }-\sigma _{j}^{\nu
}\mid
\end{equation*}
and note $\xig_{j}^{\nu }=\left( \xi _{j,1}^{\nu },\xi _{j,2}^{\nu },\xi
_{j,3}^{\nu }\right) _{\Re _{j}^{\nu }}.$
We also note $\mathbb{B}(\sigma _{j}^{\nu })$ the curvature
matrix of $\Gamma_j $ in $\sigma _{j}^{\nu }$ : it's the diagonal
matrix with diagonal entries $\left( \frac{1}{U_{j}^{\nu }},\frac{1}{V_{j}^{\nu }}%
,0\right) $ in $\Re _{j}^\nu$
where $U_j^\nu=1/k_1^j(\sigma _{j}^{\nu })$ and  $V_j^\nu=1/k_2^j(\sigma _{j}^{\nu })$
are the principal radius of curvature. Due to strict convexity of compact we have 
$U_j^\nu>0$ and $V_j^\nu>0.$

\item  Let $\mathbb{R}_{j}^{\nu }$ defined by :
\begin{equation*}
\mathbb{R}_{j}^{\nu }=\begin{pmatrix}
\DOT{\ug_{j+1}^{\nu }}{\ug_{j}^{\nu }} & 
\DOT{\ug_{j+1}^{\nu }}{\vg_{j}^{\nu }} & 
\DOT{\ug_{j+1}^{\nu }}{\nog_{j}^{\nu }} \\
\DOT{\vg_{j+1}^{\nu }}{\ug_{j}^{\nu }} & 
\DOT{\vg_{j+1}^{\nu }}{\vg_{j}^{\nu }} & 
\DOT{\vg_{j+1}^{\nu }}{\nog_{j}^{\nu }} \\
\DOT{\nog_{j+1}^{\nu }}{\ug_{j}^{\nu }} & 
\DOT{\nog_{j+1}^{\nu }}{\vg_{j}^{\nu }} & 
\DOT{\nog_{j+1}^{\nu }}{\nog_{j}^{\nu }}
\end{pmatrix}
\end{equation*}

\item  we call $\mathcal{R}_{l}(x)$ the set of all $l$-uplet
$\rho =(\sigma _{1},\cdots ,\sigma _{l})\in (\partial \Omega )_{*}^{l}$ such
that $\rho $ is an optic ray going through $x$ and $\sigma _{i}$ the point of $i^{th}$
reflection along this ray.
We note $\mathcal{R}(x)=\bigcup_{l}\mathcal{R}_{l}(x)$.

\item  Let $x\in \Omega $ and $\rho =(\sigma _{1},\cdots ,\sigma _{l})\in
\mathcal{C}_{l}(x)$. We say that $\rho$ realize :
\begin{enumerate}
\item  \textit{a transmission condition at point }$\sigma _{i}$ $(i=1,\cdots
,l)$ if
\begin{equation*}
\left\{
\begin{array}{lcll}
{\frac{\sigma _{2}-\sigma _{1}}{\mid \sigma _{2}-\sigma _{1}\mid }} & = & \xig
& \ \ \mbox{for }i=1 \\
{\frac{\sigma _{i+1}-\sigma _{i}}{\mid \sigma _{i+1}-\sigma _{i}\mid }} & =
& {\frac{\sigma _{i}-\sigma _{i-1}}{\mid \sigma _{i}-\sigma _{i-1}\mid }} &
\ \ \mbox{for }i\in \{2,\cdots ,l\}
\end{array}
\right.
\end{equation*}

\item  \textit{a reflection condition at point }$\sigma _{i}$ $(i=1,\cdots
,l)$ if
\begin{equation*}
\left\{
\begin{array}{lcl}
{\frac{\sigma _{2}-\sigma _{1}}{\mid \sigma _{2}-\sigma _{1}\mid }} & = &
\xig -2\DOT{\xig}{\nog(\sigma _{1})} \nog(\sigma _{1})\ \ \
\mbox{for }i=1 \\
{\frac{\sigma _{i+1}-\sigma _{i}}{\mid \sigma _{i+1}-\sigma _{i}\mid }} & =
& {\frac{\sigma _{i}-\sigma _{i-1}}{\mid \sigma _{i}-\sigma _{i-1}\mid }}%
-2\DOT{\frac{\sigma _{i}-\sigma _{i-1}}{\mid \sigma _{i}-\sigma_{i-1}\mid }}
{\nog(\sigma _{i})}\nog(\sigma _{i})\ \ \mbox{for }i\in
\{2,\cdots ,l\}
\end{array}
\right.
\end{equation*}
\end{enumerate}
Here $\sigma_{l+1}=x.$

\item  Let $x\in \Omega $ and $\nu =(\sigma _{1}^{\nu },\cdots ,\sigma
_{l}^{\nu })\in \mathcal{C}_{l}(x)$ we note
\begin{equation*}
\delta ^{\nu }(\sigma _{j}^{\nu })=\left\{
\begin{array}{l}
0\ \mbox{if }\sigma _{j}^{\nu }\mbox{ is a transmission point} \\
1\ \mbox{if }\sigma _{j}^{\nu }\mbox{ is a reflexion point}
\end{array}
\right.
\end{equation*}

\item  We note $\mathcal{T}=\bigcup_{l}\mathcal{T}_{l}$ with $\mathcal{T}%
_{l} $ the set of points $x\in \Omega $ such that exists $(\sigma
_{1},\cdots ,\sigma _{l})\in (\partial \Omega )_{*}^{l}$ verifying

\begin{enumerate}
\item  $\DOT{ \xig}{\nog(\sigma _{1})} =0$\ \textbf{or}\ $%
\exists j\in \{1,\cdots ,l-1\}$ such that $\DOT{\sigma _{j+1}-\sigma
_{j}}{\nog(\sigma _{j+1})}=0$,

\item  The phase $\psi _{l}(x,\bullet )$ is stationary on $(\partial \Omega
)_{*}^{l}$ in $(\sigma _{1},\cdots ,\sigma _{l})$.
\end{enumerate}
\end{itemize}

\subsection{Matrix applications}

\begin{itemize}
\item  Let $\sigma >0$, we note $\ES_{\sigma }\subset \mathcal{M}%
_{3,3}(\mathbb{R})$ the set of matrix $\mathbb{A}$ such that $\mathbb{I}%
+\sigma \mathbb{A}$ is regular. We note $S_{\sigma }$ the following
application :
\begin{equation*}
\begin{array}{ccccl}
S_{\sigma } & : & \ES_{\sigma } & \longrightarrow & \mathcal{M}%
_{3,3}(\mathbb{R}) \\
&  & \mathbb{A} & \longmapsto & \mathbb{A}(\mathbb{I}+\sigma \mathbb{A})^{-1}
\end{array}
\end{equation*}

\item  Let $\mathbb{B}\in \mathcal{M}_{3,3}(\mathbb{R})$ a symmetric matrix
, $\eg \in \mathbb{R}^{3}$, $\zg \in \mathbb{R}^{3}$. We suppose $%
\DOT{\zg}{\eta} \neq 0$. We note $T_{\mathbb{B},\eg
,\zg }$ the application of $\mathcal{M}_{3,3}(\mathbb{R})$ given
by :\newline
$\forall \mathbb{A}\in \mathcal{M}_{3,3}(\mathbb{R}),\ \forall x\in \mathbb{R}%
^{3}$
\begin{equation*}
\begin{array}{ccl}
(T_{\mathbb{B},\eg ,\zg }(\mathbb{A}))x & = & (\mathbb{A}-2
\DOT{\zg}{\eg} \mathbb{B})x-2\DOT{ \eg}{x} (
\mathbb{A}\eg +\mathbb{B}\zg ) \\
&  & -2\DOT{\mathbb{A}\eg +\mathbb{B}\zg}{x}\eg 
     +2\left[ 2\DOT{\mathbb{A}\eg}{\eg} 
              -{\frac{\DOT{\mathbb{B}\zg}{\zg }}
                     {\DOT{ \zg}{\eg} } }
       \right] \DOT{\eg}{x} \eg
\end{array}
\end{equation*}

\item  we define, for $x\in \Omega \setminus \mathcal{T}$ and $\nu
=(\sigma _{1}^{\nu },\cdots ,\sigma _{l}^{\nu })\in \mathcal{C}_{l}(x)$, the
following $l$ matrices $\mathbb{M}_{j}^{\nu }$ in $\mathcal{M}_{2,2}(\mathbb{%
R})$:
\begin{equation*}
\mathbb{M}_{1}^{\nu }=\EHb_{1,1}\psi_{l}(x,\nu )
\end{equation*}
and, $\ \forall j\in \{2,\cdots ,l\}$%
\begin{equation*}
\mathbb{M}_{j}^{\nu }=\EHb_{j,j}\psi_{l}(x;\nu )
-\EHb_{j-1,j}\psi _{l}(x;\nu )%
\left[ \mathbb{M}_{j-1}^{\nu }\right] ^{-1}\EHb_{j,j-1}\psi _{l}(x;\nu ).
\end{equation*}

\begin{remark}
We shall see in Lemma \ref{Lemma3} that $\mathbb{M}_{j}^{\nu }$ is regular.
\end{remark}

\item  Let $x\in \Omega \setminus \mathcal{T}$ and $\nu =(\sigma _{1}^{\nu
},\cdots ,\sigma _{l}^{\nu })\in \mathcal{C}_{l}(x)$. We define by
recurrence the $l$ symmetric matrices $\mathbb{P}_{j}^{\nu }$ in 
$\EL(\R^3)$
such that
\begin{equation*}
\mathbb{P}_{1}^{\nu }=T_{\mathbb{B}(\sigma _{1}^{\nu }),n(\sigma _{1}^{\nu
}),\xi }\left( 0\right) \times \delta ^{\nu }(\sigma _{1}^{\nu })
\end{equation*}
and, $\forall j\in \{2,\cdots ,l\}$%
\begin{equation*}
\mathbb{P}_{j}^{\nu }=S_{\lambda _{j-1}^{\nu }}(\mathbb{P}_{j-1}^{\nu
})\times \left( 1-\delta ^{\nu }(\sigma _{j}^{\nu })\right) +T_{\mathbb{B}%
(\sigma _{j}^{\nu }),n(\sigma _{j}^{\nu }),\xi _{j-1}^{\nu }}\left(
S_{\lambda _{j-1}^{\nu }}(\mathbb{P}_{j-1}^{\nu })\right) \times \delta
^{\nu }(\sigma _{j}^{\nu }).
\end{equation*}

\item  Let $\mathcal{D}_{l}$ the function define on $\Gamma_*^l$ by
\begin{equation*}
\begin{array}{c}
\mathcal{D}_{l}(\sigma _{1},\cdots ,\sigma _{l}) \\
= \\
\left( \left| \DOT{\xig}{\nog(\sigma _{1})} \right|
       -\DOT{\xig}{\nog(\sigma _{1})} 
\right)
\prod\limits_{j=1}^{l-1}
\frac{
  \left[ \left| 
    \DOT{ \frac{\sigma_{j+1}-\sigma _{j} }{\mid \sigma _{j+1}-\sigma _{j}\mid } }
        {\nog(\sigma_{j+1})} \right|
    -\DOT{\frac{\sigma _{j+1}-\sigma _{j}}{\mid \sigma _{j+1}-\sigma _{j}\mid } }
         {\nog(\sigma _{j+1})} 
  \right]
}{\mid \sigma _{j+1}-\sigma _{j}\mid }
\end{array}
\end{equation*}

\item  Let $\mathcal{N}_{l}$ the function define on $\Gamma_*^l$ by
\begin{equation*}
\begin{array}{c}
\mathcal{N}_{l}(\sigma _{1},\cdots ,\sigma _{l}) \\
= \\
\left( 
  \frac{ \DOT{\xig}{\nog(\sigma _{1})} }
       {\left| \DOT{\xig}{\nog(\sigma _{1})} \right| }
  -1\right)
\prod\limits_{j=1}^{l-1}
\left[ 
  \frac{ \DOT{ {\frac{\sigma_{j+1}-\sigma _{j}}{\mid \sigma _{j+1}-\sigma _{j}\mid }}}
             {\nog(\sigma_{j+1}) } }
       {\left| \DOT{ {\frac{\sigma_{j+1}-\sigma _{j}}{\mid \sigma _{j+1}-\sigma _{j}\mid }}}
             {\nog(\sigma_{j+1}) }  \right| }%
-1\right] 
\DOT{{\frac{\sigma _{j+1}-\sigma _{j}}{\mid \sigma
_{j+1}-\sigma _{j}\mid }}}{\nog(\sigma _{j})}
\end{array}
\end{equation*}
\end{itemize}

\section{Geometrical optics approximation}

We only give the main results. The geometrical optic approximation is given,
for Dirichlet and Neumann problems, respectively by : $\forall x\in \Omega
\setminus \mathcal{T}$
\begin{equation}
\begin{array}{c}
v_{O.G.}^{D}(x) \\
= \\
e^{-ik\DOT{\xig}{x} }+\sum\limits_{l \geq 1}
(-1)^{l}\sum\limits_{\rho =(\sigma _{1}^{\rho },\cdots ,\sigma _{l}^{\rho
})\in \mathcal{R}_{l}(x)}\frac{e^{-ik\psi _{l}(x;\rho )}}{%
\prod\limits_{j=1}^{l}\sqrt{\det (\mathbb{I}+\mid \sigma _{j+1}^{\rho
}-\sigma _{j}^{\rho }\mid \mathbb{P}_{j}^{\rho })}}
\end{array}
\label{GOA for Dirichlet}
\end{equation}
and
\begin{equation}
\begin{array}{c}
v_{O.G.}^{N}(x) \\
= \\
e^{-ik\DOT{\xig}{x} }+\sum\limits_{l\geq 1 }
\sum\limits_{\rho =(\sigma _{1}^{\rho },\cdots ,\sigma _{l}^{\rho })\in
\mathcal{R}_{l}(x)}\frac{e^{-ik\psi _{l}(x;\rho )}}{\prod\limits_{j=1}^{l}%
\sqrt{\det (\mathbb{I}+\mid \sigma _{j+1}^{\rho }-\sigma _{j}^{\rho }\mid
\mathbb{P}_{j}^{\rho })}}
\end{array}
\label{GOA for Neumann}
\end{equation}
Here $\sigma _{l+1}=x$, $\det \left( I+\mid \sigma _{j+1}^{\rho }-\sigma
_{j}^{\rho }\mid P_{j}^{\rho }\right) $ is positive, and in corollary, the
Maslov indice vanished. For more explanation, report to \cite{Balabane:Bardos} or \cite{Cuv:GO:2013}.

\section{Iterative Kirchhoff approximation}

\subsection{Dirichlet problem}

We introduce the following kernels series
\begin{equation*}
p_{1}^{D}(\sigma )=ik\left( \left| \left\langle \xi ,n(\sigma )\right\rangle
\right| -\left\langle \xi ,n(\sigma )\right\rangle \right)
e^{-ik\left\langle \xi ,\sigma \right\rangle }\ \ \forall \sigma \in
\partial \Omega
\end{equation*}
and, for $\sigma \in \partial K_{j}$%
\begin{equation*}
p_{l}^{D}(\sigma )=\frac{ik}{4\pi }\int_{\partial \Omega \setminus \partial
K_{j}}p_{l-1}^{D}(\sigma ^{\prime })\left[ \left| {\frac{\sigma -\sigma
^{\prime }}{\mid \sigma -\sigma ^{\prime }\mid }}.n(\sigma )\right| -{\frac{%
\sigma -\sigma ^{\prime }}{\mid \sigma -\sigma ^{\prime }\mid }}.n(\sigma )%
\right] {\frac{e^{-ik\mid \sigma ^{\prime }-\sigma \mid }}{\mid \sigma
^{\prime }-\sigma \mid }}d\sigma ^{\prime }
\end{equation*}
We set the \textbf{iterative Kirchhoff approximation} for Dirichlet problem $%
(D)$ by
\begin{equation}
\left\{
\begin{array}{lcl}
u_{0}^{D}(x) & = & 0 \\
u_{l}^{D}(x) & = & u_{l-1}^{D}(x)+\frac{1}{4\pi }\int_{\partial \Omega
}p_{l}^{D}(\sigma ){\frac{e^{-ik\mid x-\sigma \mid }}{\mid x-\sigma \mid }}%
d\sigma
\end{array}
\right.  \label{eqD1}
\end{equation}
That is to say with previous notations

\begin{equation}
u_{l}^{D}(x)=u_{l-1}^{D}(x)+\left( {\frac{ik}{4\pi }}\right)
^{l}\int_{(\partial \Omega )_{\ast }^{l}}\mathcal{D}_{l}(\sigma _{1},\cdots
,\sigma _{l}){\frac{e^{-ik\psi _{l}(x;\sigma _{1},\cdots ,\sigma _{l})}}{%
\mid x-\sigma _{l}\mid }}d\sigma _{1}\cdots d\sigma _{l}  \label{eqD2}
\end{equation}

We state the main result comparing the iterative method describe in (\ref
{eqD1}) and the geometrical optic approximation given in (\ref{GOA for
Dirichlet}) for problem $(D)$ :

\begin{theorem}
\label{Theorem1} Let $\Omega $ an open in $\mathbb{R}^{3}$, exterior of a
regular domain $\Omega ^{\prime }$ finite and disjointed reunion of strictly
convex compacts. Let $x\in \Omega \setminus \mathcal{T}$. If $\mathcal{C}%
_{l}(x)=\emptyset $ for $l>n$ then
\begin{equation}
e^{-ik\left\langle \xi ,x\right\rangle
}+u_{n}^{D}(x)-v_{0.G.}^{D}(x)=O\left( \frac{1}{k}\right)
\end{equation}
locally uniformly in $x$.
\end{theorem}

\subsection{Neumann problem}

We introduce the following kernels series
\begin{equation*}
p_{1}^{N}(\sigma )=ik\left( \frac{\left\langle \xi ,n(\sigma )\right\rangle
}{\left| \left\langle \xi ,n(\sigma )\right\rangle \right| }-1\right)
e^{-ik\left\langle \xi ,\sigma \right\rangle }\ \ \forall \sigma \in
\partial \Omega
\end{equation*}
and, for $\sigma \in \partial K_{j}$%
\begin{eqnarray*}
p_{l}^{N}(\sigma ) &=&\frac{ik}{4\pi }\int_{\partial \Omega \setminus
\partial K_{j}}p_{l-1}^{N}(\sigma ^{\prime })\left[ \frac{\left\langle {%
\frac{\sigma -\sigma ^{\prime }}{\mid \sigma -\sigma ^{\prime }\mid }}%
,n(\sigma )\right\rangle }{\left| \left\langle {\frac{\sigma -\sigma
^{\prime }}{\mid \sigma -\sigma ^{\prime }\mid }},n(\sigma )\right\rangle
\right| }-1\right] \\
&&\times \left\langle {\frac{\sigma -\sigma ^{\prime }}{\mid \sigma -\sigma
^{\prime }\mid }},n(\sigma ^{\prime })\right\rangle {\frac{e^{-ik\mid \sigma
^{\prime }-\sigma \mid }}{\mid \sigma ^{\prime }-\sigma \mid }}d\sigma
^{\prime }
\end{eqnarray*}
We set the \textbf{iterative Kirchhoff approximation} for Neumann problem $%
(N)$ by
\begin{equation}
\left\{
\begin{array}{lcl}
u_{0}^{N}(x) & = & 0 \\
u_{l}^{N}(x) & = & u_{l-1}^{N}(x)+\frac{1}{4\pi }\int_{\partial \Omega
}p_{l}^{N}(\sigma )\left\langle {\frac{x-\sigma }{\mid x-\sigma \mid }}%
,n(\sigma )\right\rangle {\frac{e^{-ik\mid x-\sigma \mid }}{\mid x-\sigma
\mid }}d\sigma
\end{array}
\right.  \label{eqN1}
\end{equation}
That is to say with previous notations
\begin{equation}
\begin{array}{c}
u_{l}^{N}(x)-u_{l-1}^{N}(x) \\
= \\
\left( {\frac{ik}{4\pi }}\right) ^{l}\int_{(\partial \Omega )_{\ast }^{l}}%
\mathcal{N}_{l}(\sigma _{1},\cdots ,\sigma _{l})\left\langle {\frac{x-\sigma
_{l}}{\mid x-\sigma _{l}\mid }},n(\sigma _{l})\right\rangle {\frac{%
e^{-ik\psi _{l}(x;\sigma _{1},\cdots ,\sigma _{l})}}{\mid x-\sigma _{l}\mid }%
}d\sigma _{1}\cdots d\sigma _{l}
\end{array}
\label{eqN2}
\end{equation}

We state the main result comparing the iterative method describe in (\ref
{eqN1}) and the geometrical optic approximation given in (\ref{GOA for
Neumann}) for problem $(N)$ :

\begin{theorem}
\label{Theorem2} Let $\Omega $ an open in $\mathbb{R}^{3}$, exterior of a
regular domain $\Omega ^{\prime }$ finite and disjointed reunion of strictly
convex compacts. Let $x\in \Omega \setminus \mathcal{T}$. If $\mathcal{C}%
_{l}(x)=\emptyset $ for $l>n$ then
\begin{equation}
e^{-ik\left\langle \xi ,x\right\rangle
}+u_{n}^{N}(x)-v_{0.G.}^{N}(x)=O\left( \frac{1}{k}\right)
\end{equation}
locally uniformly in $x$.
\end{theorem}

\section{Technical Lemmas and properties}
To prove previous theorems we need some technical lemmas and properties.
\subsection{Stationary phase points of $\protect\psi _{l}(x;\bullet )$}

We first remark that

\begin{remark}
\label{Remark1}Let $\nu =(\sigma _{1}^{\nu },\cdots ,\sigma _{l}^{\nu })\in
(\partial \Omega )_{\ast }^{l},$ if
\begin{equation*}
\left\langle \xi ,n(\sigma _{1}^{\nu })\right\rangle \geq 0
\end{equation*}
\underline{or} if exists $j\in \{1,\cdots ,l-1\}$ such that
\begin{equation*}
\left\langle {\frac{\sigma _{j+1}^{\nu }-\sigma _{j}^{\nu }}{\mid \sigma
_{j+1}^{\nu }-\sigma _{j}^{\nu }\mid }},n(\sigma _{j+1}^{\nu })\right\rangle
\geq 0
\end{equation*}
then
\begin{equation*}
\mathcal{D}_{l}(\sigma _{1}^{\nu },\cdots ,\sigma _{l}^{\nu })=\mathcal{N}%
_{l}(\sigma _{1}^{\nu },\cdots ,\sigma _{l}^{\nu })=0
\end{equation*}
\end{remark}

To find stationary phase points on $(\partial \Omega )_{\ast }^{l}$, we have
to compute, for all $x\in \mathbb{R}^{3}$, the set of points $\nu =(\sigma
_{1}^{\nu },\cdots ,\sigma _{l}^{\nu })\in (\partial \Omega )_{\ast }^{l}$
which satisfies
\begin{equation}
\nabla _{\nu }^{S}\psi _{l}(x;\sigma _{1}^{\nu },\cdots ,\sigma _{l}^{\nu
})=0  \label{eq3}
\end{equation}
We have
\begin{equation*}
\nabla _{\nu }{\psi _{l}}(x;\sigma _{1}^{\nu },\cdots ,\sigma _{l}^{\nu
})=\left(
\begin{array}{c}
\xi -{\frac{\sigma _{2}^{\nu }-\sigma _{1}^{\nu }}{\left| \sigma _{2}-\sigma
_{1}\right| }} \\
{\frac{\sigma _{2}^{\nu }-\sigma _{1}^{\nu }}{\left| \sigma _{2}^{\nu
}-\sigma _{1}^{\nu }\right| }}-{\frac{\sigma _{3}^{\nu }-\sigma _{2}^{\nu }}{%
\left| \sigma _{3}^{\nu }-\sigma _{2}^{\nu }\right| }} \\
\vdots  \\
{\frac{\sigma _{l}^{\nu }-\sigma _{l-1}^{\nu }}{\left| \sigma _{l}^{\nu
}-\sigma _{l-1}^{\nu }\right| }}-{\frac{x-\sigma _{l}^{\nu }}{\left|
x-\sigma _{l}^{\nu }\right| }}
\end{array}
\right)
\end{equation*}
The condition (\ref{eq3}) is equivalent to the existence of $(\mu
_{1},\cdots ,\mu _{l})\in \mathbb{R}^{l}$ such that
\begin{equation*}
\nabla _{\nu }{\psi _{l}}(x;\sigma _{1}^{\nu },\cdots ,\sigma _{l}^{\nu
})=\left(
\begin{array}{c}
\mu _{1}n(\sigma _{1}^{\nu }) \\
\vdots  \\
\mu _{l}n(\sigma _{l}^{\nu })
\end{array}
\right)
\end{equation*}
By hypothesis $\mid \xi \mid =1$, so we obtain
\begin{equation*}
\left\{
\begin{array}{clcl}
& \mu _{1}=0 & \mbox{or} & \mu _{1}=2\left\langle \xi ,n(\sigma _{1}^{\nu
})\right\rangle n(\sigma _{1}^{\nu }) \\
\mbox{and}\ \ \forall j\in \{2,\cdots ,l\} &  &  &  \\
& \mu _{j}=0 & \mbox{or} & \mu _{j}=2\left\langle {\frac{\sigma _{j}^{\nu
}-\sigma _{j-1}^{\nu }}{\left| \sigma _{j}^{\nu }-\sigma _{j-1}^{\nu
}\right| }},n(\sigma _{j}^{\nu })\right\rangle n(\sigma _{j}^{\nu })
\end{array}
\right.
\end{equation*}
Then, using remark \ref{Remark1}, we have the

\begin{lemma}
\label{Lemma1} Let $x\in \Omega \setminus \mathcal{T}$ and $\nu =(\sigma
_{1}^{\nu },\cdots ,\sigma _{l}^{\nu })\in \mathcal{C}_{l}(x)$. Then $\nu $
realize a transmission or a reflection condition on each points $\sigma
_{j}^{\nu },\,j\in \{1,\ldots ,l\}.$
\end{lemma}

\subsection{Relation between $\mathbb{M}_{j}^{\protect\nu }$ and $\mathbb{P}{%
_{j}^{\protect\nu }}$}

We have the fundamental Lemma

\begin{lemma}
\label{Lemma2}Let $x\in \Omega \setminus \mathcal{T}$ and $\nu =(\sigma
_{1}^{\nu },\cdots ,\sigma _{l}^{\nu })\in C_{l}(x)$. Then, $\forall j\in
\{1,\cdots ,l\}$ we have :
\begin{equation}
\mathbb{M}_{j}^{\nu }=\mathbb{P}{_{j}^{\nu }}_{\mid T_{{\partial \Omega }}{(}%
\sigma _{j}^{\nu })}+\frac{1}{\lambda _{j}^{\nu }}\left( \mathbb{I}-(\xi
_{j}^{\nu }\ ^{t}\xi _{j}^{\nu })_{\mid T_{\partial \Omega }{(}%
\sigma _{j}^{\nu })}\right),   \label{eq4}
\end{equation}
\begin{equation}
\det \mathbb{M}_{j}^{\nu }=\left( \frac{\left\langle \xi _{j}^{\nu
},n(\sigma _{j}^{\nu })\right\rangle }{\lambda _{j}^{\nu }}\right) ^{2}\det (%
\mathbb{I}+\lambda _{j}^{\nu }\mathbb{P}_{j}^{\nu }),  \label{eq5}
\end{equation}
\begin{equation}
\sgn\mathbb{M}_{j}^{\nu }=2,  \label{eq6}
\end{equation}
and
\begin{equation}
\mathbb{M}_{j}^{\nu }\ \ \ \mbox{inversible}.  \label{eq7}
\end{equation}
\end{lemma}

\textsc{Proof of lemma }\ref{Lemma2}\textsc{\ :}\newline
The proof of the lemma worked out by recurrence.

\begin{center}
\underline{\textbf{step one of recurrence proof}}
\end{center}

Using definition of $\mathbb{P}_{1}^{\nu },$ we easily obtain
\begin{equation*}
\mathbb{P}{_{1}^{\nu }}_{\mid T_{{\partial \Omega }}{(}\sigma _{1}^{\nu
})}=2\left\langle \xi _{1}^{\nu },n(\sigma _{1}^{\nu })\right\rangle \delta
^{\nu }(\sigma _{1}^{\nu })\mathbb{B}(\sigma _{1}^{\nu })_{\mid T_{{\partial
\Omega }}{(}\sigma _{1}^{\nu })}
\end{equation*}
On the other hand, we have
\begin{equation*}
\begin{array}{lcl}
\mathbb{M}_{1}^{\nu } & = & H_{\sigma _{1},\sigma _{1}}^{S}\psi _{l}(x;\nu )
\\
& = & H_{\sigma _{1},\sigma _{1}}^{S}\left( \left\langle \xi ,\sigma
_{1}^{\nu }\right\rangle +\left| \sigma _{2}^{\nu }-\sigma _{1}^{\nu
}\right| \right) \\
& = & 2\left\langle \xi _{1}^{\nu },n(\sigma _{1}^{\nu })\right\rangle
\delta ^{\nu }(\sigma _{1}^{\nu })\mathbb{B}(\sigma _{1}^{\nu })_{\mid T_{{%
\partial \Omega }}{(}\sigma _{1}^{\nu })}+\frac{1}{\lambda _{1}^{\nu }}%
\left( \mathbb{I}-(\xi _{1}^{\nu }\ ^{t}\xi _{1}^{\nu })_{\mid T_{\partial
\Omega }(\sigma _{1}^{\nu })}\right)
\end{array}
\end{equation*}
and so we have proved formula (\ref{eq4}) for $j=1.$

To obtain formula (\ref{eq5}), we first remark that $\mathbb{P}_{1}^{\nu
}\xi _{1}^{\nu }=0$. We set
\begin{equation*}
\mathbb{H}_{j}=\left[
\begin{array}{ccc}
1 & 0 & \xi _{j,1}^{\nu } \\
0 & 1 & \xi _{j,2}^{\nu } \\
0 & 0 & \xi _{j,3}^{\nu }
\end{array}
\right]
\end{equation*}
Then, combining $\mid \xi _{1}^{\nu }\mid =1$ and formula (\ref{eq4}) gives
\begin{equation*}
\begin{array}{lcl}
\det (\mathbb{I}+\lambda _{1}^{\nu }\mathbb{P}_{1}^{\nu }) & = & \frac{1}{%
(\det \mathbb{H}_{1})^{2}}\det \left( \ ^{t}\mathbb{H}_{1}\left( \mathbb{I}%
+\lambda _{1}^{\nu }\mathbb{P}_{1}^{\nu }\right) \mathbb{H}_{1}\right) \\
& = & \left( \frac{\lambda _{1}^{\nu }}{\left\langle \xi _{1}^{\nu
},n(\sigma _{1}^{\nu })\right\rangle }\right) ^{2}\det \mathbb{M}_{1}^{\nu }.
\end{array}
\end{equation*}

One finds easily that $\det (\mathbb{I}+\lambda _{1}^{\nu }\mathbb{P}%
_{1}^{\nu })>0 $ and $\tr\mathbb{M}_{1}^{\nu }>0$. So, we get $%
\sgn\mathbb{M}_{1}^{\nu }=2$.

We deduce that $\det \mathbb{M}_{1}^{\nu }>0$ and so $\mathbb{M}%
_{1}^{\nu }$ is regular.

\begin{center}
\underline{\textbf{Step $j+1$ of recurrence proof}}
\end{center}

To prove formula (\ref{eq4}) at step $j+1$, we first have to compute $(%
\mathbb{P}_{j+1}^{\nu })_{\Re _{j+1}}$ and $\mathbb{M}_{j+1}^{\nu }$
respectively in function of $(\mathbb{P}_{j}^{\nu })_{\Re _{j}^{\nu }}$ and $%
\mathbb{M}_{j}^{\nu }. $ To simplify notations, we note $(\mathbb{P}%
_{j}^{\nu })_{\Re _{j}^{\nu }}=\left[ P_{r,s}\right] _{r,s\in \{1,2,3\}},$ and
$\mathbb{M}_{j}^{\nu }=\left[ M_{r,s}\right] _{r,s\in \{1,2\}}.$

\begin{itemize}
\item  \textbf{Computation of }$(\mathbb{P}_{j+1}^{\nu })_{\Re _{j+1}}$

By definition,
\begin{equation*}
\begin{array}{c}
(\mathbb{P}_{j+1}^{\nu })_{\Re _{j+1}^{\nu }} \\
= \\
\left( S_{\lambda _{j}^{\nu }}(\mathbb{P}_{j}^{\nu })\right) _{\Re
_{j+1}^{\nu }}\times \left( 1-\delta ^{\nu }(\sigma _{j+1}^{\nu })\right)
+T_{\mathbb{B}(\sigma _{j+1}^{\nu }),n(\sigma _{j+1}^{\nu }),\xi _{j}^{\nu
}}\left( S_{\lambda _{j}^{\nu }}(\mathbb{P}_{j}^{\nu })\right) _{\Re
_{j+1}^{\nu }}\times \delta ^{\nu }(\sigma _{j+1}^{\nu }).
\end{array}
\end{equation*}
We first evaluate $S_{\lambda _{j}^{\nu }}(\mathbb{P}_{j}^{\nu })_{\Re _{j}} :$
\begin{equation*}
(S_{\lambda _{j}^{\nu }}(\mathbb{P}_{j}^{\nu }))_{\Re _{j}}\det (\mathbb{I}%
+\lambda _{j}^{\nu }\mathbb{P}_{j}^{\nu })=(\mathbb{P}_{j}^{\nu })_{\Re
_{j}}+\lambda _{j}^{\nu }\mathbb{Q}_{j}^{\nu }+\left( \lambda _{j}^{\nu
}\right) ^{2}\mathbb{I}\det \mathbb{P}_{j}^{\nu }
\end{equation*}
with $\mathbb{Q}_{j}^{\nu }=\left[ Q_{r,s}\right] _{p,q\in \{1,2,3\}}\in
\mathcal{M}_{3,3}(\mathbb{R})$ and
\begin{equation*}
\begin{array}{lll}
Q_{11} & = & P_{11}\left( P_{33}+P_{22}\right) -P_{12}^{2}-P_{13}^{2}, \\
Q_{12} & = & P_{12}P_{33}-P_{13}P_{23}, \\
Q_{13} & = & P_{13}P_{22}-P_{12}P_{23}, \\
Q_{22} & = & P_{22}\left( P_{33}+P_{11}\right) -P_{12}^{2}-P_{23}^{2}, \\
Q_{23} & = & P_{23}P_{11}-P_{12}P_{13}, \\
Q_{33} & = & P_{33}\left( P_{22}+P_{11}\right) -P_{13}^{2}-P_{23}^{2}.
\end{array}
\end{equation*}
As $\det \mathbb{P}_{j}^{\nu }=0$, we obtain :
\begin{equation*}
(S_{\lambda _{j}^{\nu }}(\mathbb{P}_{j}^{\nu }))_{\Re _{j}^{\nu }}=\frac{1}{%
\det (\mathbb{I}+\lambda _{j}^{\nu }\mathbb{P}_{j}^{\nu })}\left( (\mathbb{P}%
_{j}^{\nu })_{\Re _{j}^{\nu }}+\lambda _{j}^{\nu }\mathbb{Q}_{j}^{\nu
}\right)
\end{equation*}
and
\begin{equation*}
\begin{array}{lcl}
(S_{\lambda _{j}^{\nu }}(\mathbb{P}_{j}^{\nu }))_{\Re _{j+1}^{\nu }} & = &
\frac{1}{\det (\mathbb{I}+\lambda _{j}^{\nu }\mathbb{P}_{j}^{\nu })}{\ }%
\mathbb{R}_{j}^{\nu }\left( (\mathbb{P}_{j}^{\nu })_{\Re _{j}}+\lambda
_{j}^{\nu }\mathbb{Q}_{j}^{\nu }\right) \,^{t}\mathbb{R}_{j}^{\nu }.
\end{array}
\end{equation*}
Now, we have to compute $T_{\mathbb{B}(\sigma _{j+1}^{\nu }),n(\sigma
_{j+1}^{\nu }),\xi _{j}^{\nu }}\left( S_{\lambda _{j}^{\nu }}(\mathbb{P}%
_{j}^{\nu })\right) _{\Re _{j+1}^{\nu }}$.

In local coordinates, we have $\mathbb{B}(\sigma _{j+1}^{\nu })=\left[
\begin{array}{lll}
\frac{1}{U_{j+1}^{\nu }} & 0 & 0 \\
0 & \frac{1}{V_{j+1}^{\nu }} & 0 \\
0 & 0 & 0
\end{array}
\right] _{\Re _{j+1}^{\nu }},$ $n(\sigma _{j+1}^{\nu })_{\Re _{j+1}^{\nu
}}=\left(
\begin{array}{l}
0 \\
0 \\
1
\end{array}
\right) $ and $\left( \xi _{j+1}^{\nu }\right) _{\Re _{j+1}^{\nu }}=\left(
\begin{array}{l}
\xi _{j+1,1}^{\nu } \\
\xi _{j+1,2}^{\nu } \\
\xi _{j+1,3}^{\nu }
\end{array}
\right) .$ So, if we note $\mathbb{W=}\left[ (w_{pq})_{p,q\in \{1,2,3\}}%
\right] =(S_{\lambda _{j}^{\nu }}(\mathbb{P}_{j}^{\nu }))_{\Re _{j+1}^{\nu
}} $ and $\mathbb{T=}\left[ (t_{pq})_{p,q\in \{1,2,3\}}\right] =\left( T_{%
\mathbb{B}(\sigma _{j+1}^{\nu }),n(\sigma _{j+1}^{\nu }),\xi _{j}^{\nu
}}\left( S_{\lambda _{j}^{\nu }}(\mathbb{P}_{j}^{\nu })\right) \right) _{\Re
_{j+1}^{\nu }}$ then, $\forall x\in \mathbb{R}^{3}$%
\begin{equation*}
\begin{array}{c}
\left( T_{\mathbb{B}(\sigma _{j+1}^{\nu }),n(\sigma _{j+1}^{\nu }),\xi
_{j}^{\nu }}\left( S_{\lambda _{j}^{\nu }}(\mathbb{P}_{j}^{\nu })\right)
\right) _{\Re _{j+1}^{\nu }}x \\
= \\
\left[
\begin{array}{lll}
w_{11}-2\frac{\left\langle \xi _{j}^{\nu },n(\sigma _{j+1}^{\nu
})\right\rangle }{U_{j+1}^{\nu }} & w_{12} & w_{13} \\
w_{12} & w_{22}-2\frac{\left\langle \xi _{j}^{\nu },n(\sigma _{j+1}^{\nu
})\right\rangle }{V_{j+1}^{\nu }} & w_{23} \\
w_{13} & w_{23} & w_{33}
\end{array}
\right] \left[
\begin{array}{l}
x_{1} \\
x_{2} \\
x_{3}
\end{array}
\right] \\
-2x_{3}\left[
\begin{array}{l}
w_{13}+\frac{\left\langle \xi _{j}^{\nu },\mathbf{u}_{j+1}^{\nu
}\right\rangle }{U_{j+1}^{\nu }} \\
w_{23}+\frac{\left\langle \xi _{j}^{\nu },\mathbf{v}_{j+1}^{\nu
}\right\rangle }{V_{j+1}^{\nu }} \\
w_{33}
\end{array}
\right] \\
-2\left( \left[
\begin{array}{l}
w_{13}+\frac{\left\langle \xi _{j}^{\nu },\mathbf{u}_{j+1}^{\nu
}\right\rangle }{U_{j+1}^{\nu }} \\
w_{23}+\frac{\xi _{j,2}^{\nu }}{V_{j+1}^{\nu }} \\
w_{33}
\end{array}
\right] ,\left[
\begin{array}{l}
x_{1} \\
x_{2} \\
x_{3}
\end{array}
\right] \right) \left[
\begin{array}{l}
0 \\
0 \\
1
\end{array}
\right] \\
+4\left( \left[
\begin{array}{l}
w_{13} \\
w_{23} \\
w_{33}
\end{array}
\right] ,\left[
\begin{array}{l}
0 \\
0 \\
1
\end{array}
\right] \right) x_{3}\left[
\begin{array}{l}
0 \\
0 \\
1
\end{array}
\right] \\
-\frac{2}{\left\langle \xi _{j}^{\nu },n(\sigma _{j+1}^{\nu })\right\rangle }%
\left( \frac{\left\langle \xi _{j}^{\nu },\mathbf{u}_{j+1}^{\nu
}\right\rangle ^{2}}{U_{j+1}^{\nu }}+\frac{\left\langle \xi _{j}^{\nu },%
\mathbf{v}_{j+1}^{\nu }\right\rangle ^{2}}{V_{j+1}^{\nu }}\right) x_{3}\left[
\begin{array}{l}
0 \\
0 \\
1
\end{array}
\right]
\end{array}
\end{equation*}
thus
\begin{equation*}
\begin{array}{l}
t_{11}=w_{11}-2\frac{\left\langle \xi _{j}^{\nu },n(\sigma _{j+1}^{\nu
})\right\rangle }{U_{j+1}^{\nu }}, \\
t_{12}=w_{12}, \\
t_{13}=-w_{13}-2\frac{\left\langle \xi _{j}^{\nu },\mathbf{u}_{j+1}^{\nu
}\right\rangle }{U_{j+1}^{\nu }}, \\
t_{22}=w_{22}-2\frac{\left\langle \xi _{j}^{\nu },n(\sigma _{j+1}^{\nu
})\right\rangle }{V_{j+1}^{\nu }}, \\
t_{23}=-w_{23}-2\frac{\left\langle \xi _{j}^{\nu },\mathbf{v}_{j+1}^{\nu
}\right\rangle }{V_{j+1}^{\nu }}, \\
t_{33}=w_{33}-\frac{2}{\left\langle \xi _{j}^{\nu },n(\sigma _{j+1}^{\nu
})\right\rangle }\left( \frac{\left\langle \xi _{j}^{\nu },\mathbf{u}%
_{j+1}^{\nu }\right\rangle ^{2}}{U_{j+1}^{\nu }}+\frac{\left\langle \xi
_{j}^{\nu },\mathbf{v}_{j+1}^{\nu }\right\rangle ^{2}}{V_{j+1}^{\nu }}\right).
\end{array}
\end{equation*}
Using $\left\langle \xi _{j}^{\nu },n(\sigma _{j+1}^{\nu })\right\rangle
=\left( 1-2\delta ^{\nu }(\sigma _{j+1}^{\nu })\right) \left\langle \xi
_{j+1}^{\nu },n(\sigma _{j+1}^{\nu })\right\rangle,$ we obtain :
\begin{equation*}
\begin{array}{c}
\left( T_{\mathbb{B}(\sigma _{j+1}^{\nu }),n(\sigma _{j+1}^{\nu }),\xi
_{j}^{\nu }}\left( S_{\lambda _{j}^{\nu }}(\mathbb{P}_{j}^{\nu })\right)
\right) _{\left| T_{{\partial \Omega }}{(}\sigma _{j+1}^{\nu })\right. } \\
= \\
\left( (S_{\lambda _{j}^{\nu }}(\mathbb{P}_{j}^{\nu }))_{\Re _{j+1}}\right)
_{\left| T_{{\partial \Omega }}{(}\sigma _{j+1}^{\nu })\right. }-2\left(
1-2\delta ^{\nu }(\sigma _{j+1}^{\nu })\right) \left\langle \xi _{j+1}^{\nu
},n(\sigma _{j+1}^{\nu })\right\rangle \mathbb{B}(\sigma _{j+1}^{\nu
})_{\left| T_{{\partial \Omega }}{(}\sigma _{j+1}^{\nu })\right. }
\end{array}
\end{equation*}
Finally, we have
\begin{equation}
\begin{array}{c}
\left( (\mathbb{P}_{j+1}^{\nu })_{\Re _{j+1}}\right) _{\left| T_{{\partial
\Omega }}{(}\sigma _{1+1}^{\nu })\right. } \\
= \\
\left( (S_{\lambda _{j}^{\nu }}(\mathbb{P}_{j}^{\nu }))_{\Re _{j+1}}\right)
_{\left| T_{{\partial \Omega }}{(}\sigma _{j+1}^{\nu })\right. } \\
+2\delta ^{\nu }(\sigma _{j+1}^{\nu })\left\langle \xi _{j+1}^{\nu
},n(\sigma _{j+1}^{\nu })\right\rangle \mathbb{B}(\sigma _{j+1}^{\nu
})_{\left| T_{{\partial \Omega }}{(}\sigma _{j+1}^{\nu })\right. }
\end{array}
\label{eq8}
\end{equation}

\item  \textbf{Computation of }$\mathbb{M}_{j+1}^{\nu }$

By definition,
\begin{equation*}
\mathbb{M}_{j+1}^{\nu }=H_{\sigma _{j+1}^{\nu },\sigma _{j+1}^{\nu
}}^{S}\psi _{l}(x,\nu )-H_{\sigma _{j}^{\nu },\sigma _{j+1}^{\nu }}^{S}\psi
_{l}(x,\nu )\left[ \mathbb{M}_{j}^{\nu }\right] ^{-1}H_{\sigma _{j+1}^{\nu
},\sigma _{j}^{\nu }}^{S}\psi _{l}(x,\nu )
\end{equation*}
We first evaluate $H_{\sigma _{j+1}^{\nu },\sigma _{j}^{\nu }}^{S}\psi
_{l}(x;\nu )[\mathbb{M}_{j}^{\nu }]^{-1}H_{\sigma _{j}^{\nu },\sigma
_{j+1}^{\nu }}^{S}\psi _{l}(x;\nu )$. By construction of $\psi _{l}$, we
have :
\begin{equation*}
H_{\sigma _{j+1}^{\nu },\sigma _{j}^{\nu }}^{S}\psi _{l}(x;\nu )=H_{\sigma
_{j+1}^{\nu },\sigma _{j}^{\nu }}^{S}\mid \sigma _{j+1}^{\nu }-\sigma
_{j}^{\nu }\mid .
\end{equation*}
That's give in $\Re _{j}$
\begin{equation*}
H_{\sigma _{j+1}^{\nu },\sigma _{j}^{\nu }}^{S}\psi _{l}(x;\nu )=\frac{1}{%
\lambda _{j}^{\nu }}\mathbb{L}_{j}^{\nu }
\end{equation*}
where $\mathbb{L}_{j}^{\nu }=\left[ L_{p,q}\right] _{_{p,q\in \{1,2\}}}$
with
\begin{eqnarray*}
L_{11} &=&\left\langle \mathbf{u}_{j}^{\nu },\xi _{j}^{\nu }\right\rangle
\left\langle \mathbf{u}_{j+1}^{\nu },\xi _{j}^{\nu }\right\rangle
-\left\langle \mathbf{u}_{j}^{\nu },\,\mathbf{u}_{j+1}^{\nu }\right\rangle, \\
L_{12} &=&\left\langle \mathbf{v}_{j}^{\nu },\xi _{j}^{\nu }\right\rangle
\left\langle \mathbf{u}_{j+1}^{\nu },\xi _{j}^{\nu }\right\rangle
-\left\langle \mathbf{v}_{j}^{\nu },\,\mathbf{u}_{j+1}^{\nu }\right\rangle, \\
L_{21} &=&\left\langle \mathbf{u}_{j}^{\nu },\xi _{j}^{\nu }\right\rangle
\left\langle \mathbf{v}_{j+1}^{\nu },\xi _{j}^{\nu }\right\rangle
-\left\langle \mathbf{u}_{j}^{\nu },\mathbf{v}_{j+1}^{\nu }\right\rangle, \\
L_{22} &=&\left\langle \mathbf{v}_{j}^{\nu },\xi _{j}^{\nu }\right\rangle
\left\langle \mathbf{v}_{j+1}^{\nu },\xi _{j}^{\nu }\right\rangle
-\left\langle \mathbf{v}_{j}^{\nu },\,\mathbf{v}_{j+1}^{\nu }\right\rangle.
\end{eqnarray*}
We have by hypothesis $\mathbb{M}_{j}^{\nu }=\mathbb{P}{_{j}^{\nu }}_{\mid
T_{{\partial \Omega }}{(}\sigma _{j}^{\nu })}+\frac{1}{\lambda _{j}^{\nu }}%
\left( \mathbb{I}-(\xi _{j}^{\nu }\ ^{t}\xi _{j}^{\nu })\right) _{\mid T_{{%
\partial \Omega }}{(}\sigma _{j}^{\nu })}$ i.e. :
\begin{eqnarray*}
\left( \mathbb{M}_{j}^{\nu }\right) _{_{\Re _{j}^{\nu }}} &=&\left[ M_{p,q}%
\right] _{p,q\in \{1,2\}} \\
&=&\left[
\begin{array}{ll}
P_{11} & P_{12} \\
P_{12} & P_{22}
\end{array}
\right] +\frac{1}{\lambda _{j}^{\nu }}\left[
\begin{array}{ll}
1-\left( \xi _{j,1}^{\nu }\right) ^{2} & -\xi _{j,1}^{\nu }\xi _{j,2}^{\nu }
\\
-\xi _{j,1}^{\nu }\xi _{j,2}^{\nu } & 1-\left( \xi _{j,2}^{\nu }\right) ^{2}
\end{array}
\right]
\end{eqnarray*}
By recurrence hypothesis $\mathbb{M}_{j}^{\nu }$ is regular, so we have :
\begin{equation*}
H_{\sigma _{j+1}^{\nu },\sigma _{j}^{\nu }}^{S}\psi _{l}(x;\nu )[\mathbb{M}%
_{j}^{\nu }]^{-1}H_{\sigma _{j}^{\nu },\sigma _{j+1}^{\nu }}^{S}\psi
_{l}(x;\nu )=\frac{1}{\left( \lambda _{j}^{\nu }\right) ^{2}}\mathbb{L}%
_{j}^{\nu }[\mathbb{M}_{j}^{\nu }]^{-1}\,^{t}\mathbb{L}_{j}^{\nu }.
\end{equation*}

Now, we evaluate $H_{\sigma _{j+1}^{\nu },\sigma _{j+1}^{\nu }}^{S}\psi
_{l}(x;\nu )$ in $\Re _{j+1}.$ We find
\begin{equation*}
\begin{array}{c}
H_{\sigma _{j+1}^{\nu },\sigma _{j+1}^{\nu }}^{S}\psi _{l}(x;\nu ) \\
= \\
H_{\sigma _{j+1}^{\nu },\sigma _{j+1}^{\nu }}^{S}\left( \mid \sigma
_{j+1}^{\nu }-\sigma _{j}^{\nu }\mid +\mid \sigma _{j+2}^{\nu }-\sigma
_{j+1}^{\nu }\mid \right) \\
= \\
\frac{1}{\lambda _{j}^{\nu }}\left( \mathbb{I}-\xi _{j}^{\nu }\,^{t}\xi
_{j}^{\nu }\right) _{\left| T_{{\partial \Omega }}{(}\sigma _{j+1}^{\nu
})\right. }+\frac{1}{\lambda _{j+1}^{\nu }}\left( \mathbb{I}-\xi _{j+1}^{\nu
}\,^{t}\xi _{j+1}^{\nu }\right) _{\left| T_{{\partial \Omega }}{(}\sigma
_{j+1}^{\nu })\right. } \\
-2\delta ^{\nu }(\sigma _{j+1}^{\nu })\left\langle n(\sigma _{j+1}^{\nu
}),\xi _{j}^{\nu }\right\rangle \mathbb{B}(\sigma _{j+1}^{\nu })_{\left| T_{{%
\partial \Omega }}{(}\sigma _{j+1}^{\nu })\right. }.
\end{array}
\end{equation*}
Finally, we obtain
\begin{equation}
\begin{array}{c}
\mathbb{M}_{j+1}^{\nu } \\
= \\
\frac{1}{\lambda _{j}^{\nu }}\left( \mathbb{I}-\xi _{j}^{\nu }\,^{t}\xi
_{j}^{\nu }\right) _{\left| T_{{\partial \Omega }}{(}\sigma _{j+1}^{\nu
})\right. }+\frac{1}{\lambda _{j+1}^{\nu }}\left( \mathbb{I}-\xi _{j+1}^{\nu
}\,^{t}\xi _{j+1}^{\nu }\right) _{\left| T_{{\partial \Omega }}{(}\sigma
_{j+1}^{\nu })\right. } \\
-2\delta ^{\nu }(\sigma _{j+1})\left\langle n(\sigma _{j+1}^{\nu }),\xi
_{j}^{\nu }\right\rangle \mathbb{B}(\sigma _{j+1})_{\left| T_{{\partial
\Omega }}{(}\sigma _{j+1}^{\nu })\right. }-\frac{1}{\left( \lambda _{j}^{\nu
}\right) ^{2}}\mathbb{L}_{j}^{\nu }[\mathbb{M}_{j}^{\nu }]^{-1}\,^{t}\mathbb{%
L}_{j}^{\nu }.
\end{array}
\label{eq9}
\end{equation}

\item  \textbf{Formula (\ref{eq4}) at step $j+1$ : } \newline
We compute now $\mathbb{D}{_{j+1}^{\nu }}=\mathbb{P}{_{j+1}^{\nu }}_{\left|
T_{{\partial \Omega }}{(}\sigma _{j+1}^{\nu })\right. }+\frac{1}{\lambda
_{j+1}^{\nu }}\left( \mathbb{I}-(\xi _{j+1}\ ^{t}\xi _{j+1})_{\left| T_{{%
\partial \Omega }}{(}\sigma _{j+1}^{\nu })\right. }\right) -\mathbb{M}%
_{j+1}^{\nu }.$ Combining formulas (\ref{eq8}) et (\ref{eq9}), we obtain
\begin{equation*}
\begin{array}{c}
\mathbb{D}{_{j+1}^{\nu }} \\
= \\
\left( (S_{\lambda _{j}^{\nu }}(\mathbb{P}_{j}^{\nu }))_{\Re _{j+1}}\right)
_{\left| T_{{\partial \Omega }}{(}\sigma _{j+1}^{\nu })\right. }+\frac{1}{%
\lambda _{j}^{\nu }}\left( \mathbb{I}-\xi _{j}^{\nu }\,^{t}\xi _{j}^{\nu
}\right) _{\left| T_{{\partial \Omega }}{(}\sigma _{j+1}^{\nu })\right. }+%
\frac{1}{\left( \lambda _{j}^{\nu }\right) ^{2}}\mathbb{L}_{j}^{\nu }[%
\mathbb{M}_{j}^{\nu }]^{-1}\,^{t}\mathbb{L}_{j}^{\nu }.
\end{array}
\end{equation*}
We prove that $\mathbb{D}{_{j+1}^{\nu }}=0$ using relations $\mathbb{P}%
_{j}^{\nu }\xi _{j}^{\nu }=0$, $\mid \xi _{j}^{\nu }\mid =1$ and $\ {^{t}}%
\mathbb{R}_{j}^{\nu }\mathbb{R}_{j}^{\nu }=\mathbb{I}$.

\item \textbf{Formula (\ref{eq5}) at step $j+1$ : } \newline
We remark that $\mathbb{P}_{j+1}^{\nu }\xi _{j+1}^{\nu }=0$. Then, we pose
\begin{equation*}
\mathbb{H}_{j+1}^{\nu }=\left[
\begin{array}{ccc}
1 & 0 & \xi _{j+1,1}^{\nu } \\
0 & 1 & \xi _{j+1,2}^{\nu } \\
0 & 0 & \xi _{j+1,3}^{\nu }
\end{array}
\right]
\end{equation*}
and we obtain $$\det \left( \mathbb{I}+\lambda _{j+1}^{\nu }\mathbb{P}%
_{j+1}^{\nu }\right) =\frac{1}{(\det \mathbb{H}_{j+1}^{\nu })^{2}}\det
\left( \ ^{t}\mathbb{H}_{j+1}^{\nu }\left( \mathbb{I}+\lambda _{j+1}^{\nu }%
\mathbb{P}_{j+1}^{\nu }\right) \mathbb{H}_{j+1}^{\nu }\right).$$
So
\begin{equation*}
\begin{array}{c}
\left\langle \xi _{{j+1}}^{\nu },n(\sigma _{j+1}^{\nu })\right\rangle
^{2}\det \left( \mathbb{I}+\lambda _{j+1}^{\nu }\mathbb{P}_{j+1}^{\nu
}\right) \\
= \\
\det \left[
\begin{array}{ccc}
1+\lambda _{j+1}^{\nu }\left\langle \mathbb{P}_{j+1}^{\nu }\mathbf{u}%
_{j+1}^{\nu },\mathbf{u}_{j+1}^{\nu }\right\rangle & \lambda _{j+1}^{\nu
}\left\langle \mathbb{P}_{j+1}^{\nu }\mathbf{u}_{j+1}^{\nu },\mathbf{v}%
_{j+1}^{\nu }\right\rangle & \left\langle \xi _{{j+1}}^{\nu },\mathbf{u}%
_{j+1}^{\nu }\right\rangle \\
\lambda _{j+1}^{\nu }\left\langle \mathbb{P}_{j+1}^{\nu }\mathbf{u}%
_{j+1}^{\nu },\mathbf{v}_{j+1}^{\nu }\right\rangle & 1+\lambda _{j+1}^{\nu
}\left\langle \mathbb{P}_{j+1}^{\nu }\mathbf{v}_{j+1}^{\nu },\mathbf{v}%
_{j+1}^{\nu }\right\rangle & \left\langle \xi _{{j+1}}^{\nu },\mathbf{v}%
_{j+1}^{\nu }\right\rangle \\
\left\langle \xi _{{j+1}}^{\nu },\mathbf{u}_{j+1}^{\nu }\right\rangle &
\left\langle \xi _{{j+1}}^{\nu },\mathbf{v}_{j+1}^{\nu }\right\rangle & 1
\end{array}
\right] .
\end{array}
\end{equation*}
As $\mid \xi _{j+1}^{\nu }\mid =1$, we get, with formula (\ref{eq4}) at step
${j+1}$:
\begin{equation*}
\det \left( \mathbb{I}+\lambda _{j+1}^{\nu }\mathbb{P}_{j+1}^{\nu }\right)
=\left( \frac{\lambda _{j+1}^{\nu }}{\left\langle \xi _{{j+1}}^{\nu
},n(\sigma _{j+1}^{\nu })\right\rangle }\right) ^{2}\det \mathbb{M}%
_{j+1}^{\nu }.
\end{equation*}
\item \textbf{Formulas (\ref{eq6}) and (\ref{eq7}) at step $j+1$ : } \newline

In fact, we proved that
\begin{equation*}
\forall \lambda >0\ \ \det (\mathbb{I}+\lambda \mathbb{P}_{j+1}^{\rho
})=\left( \frac{\lambda }{\left\langle \xi _{{j+1}}^{\rho },n(\sigma
_{j+1}^{\rho })\right\rangle }\right) ^{2}\det \mathbb{M}_{j+1}^{\rho }
\end{equation*}
where $\rho =(\sigma _{1}^{\nu },\cdots ,\sigma _{j+1}^{\nu })\in \mathcal{C}%
_{j+1}(\lambda \xi _{j+1}^{\nu }).$
Moreover, we have :
$$\mathbb{P}_{j+1}^{\rho }=\mathbb{P}_{j+1}^{\nu },$$ for $\lambda
=\lambda _{j+1}^{\nu }$, $$\mathbb{M}_{j+1}^{\rho }=\mathbb{M}_{j+1}^{\nu }$$
and $\forall \lambda >0$
$$\det (\mathbb{I}+\lambda \mathbb{P}%
_{j+1}^{\rho })>0$$ because $\mathbb{P}_{j+1}^{\rho }$ is a positive matrix.
That's give
\begin{equation*}
\forall \lambda >0,\ \ \det \mathbb{M}_{j+1}^{\rho }\geq 0.
\end{equation*}
This quantity is positive for $\lambda $ in a neighborhood of zero, then by
continuity we have :
\begin{equation*}
\forall \lambda >0\ \ \tr \mathbb{M}_{j+1}^{\rho }>0
\end{equation*}
with
\begin{equation*}
\tr \mathbb{M}_{j+1}^{\rho }=P_{11}+P_{22}+\frac{1+\left\langle \xi
_{{j+1}}^{\nu },n(\sigma _{j+1}^{\nu })\right\rangle ^{2}}{\lambda }.
\end{equation*}
Thus, we obtain formula (\ref{eq6}) to step $j+1$.\newline
We have also showed that $\mathbb{M}_{j+1}^{\nu }$ is regular.
\end{itemize}
That's close the proof of Lemma \ref{Lemma2}$\Box $

\subsection{Lemma of transmission}

\begin{lemma}[of transmission]
\label{Lemma3} 
Let $x\in \Omega \setminus \mathcal{T}$ and $\nu =(\sigma _{1}^{\nu },\cdots
,\sigma _{l}^{\nu })\in \mathcal{C}_{l}(x)$

\begin{enumerate}
\item  if exists $t\in \mathbb{R}_{\ast }^{+}$ such that
\begin{equation*}
\sigma =\sigma _{1}^{\nu }-t\xi \ \in \partial \Omega \ \mbox{and }\xi
.n(\sigma )<0
\end{equation*}
then
\begin{equation*}
\mu =(\sigma ,\sigma _{1}^{\nu },\cdots ,\sigma _{l}^{\nu })\in \mathcal{C}%
_{l+1}(x)
\end{equation*}

\item  if exists $j\in \{1,\cdots ,l-1\}$ such that
\begin{equation*}
\sigma \in \left\{ ]\sigma _{j}^{\nu };\sigma _{j+1}^{\nu }[\cap \partial
\Omega \right\} \ \mbox{and }\ (\sigma -\sigma _{j}^{\nu }).n(\sigma )<0
\end{equation*}
then
\begin{equation*}
\mu =(\sigma _{1}^{\nu },\cdots ,\sigma _{j}^{\nu },\sigma ,\sigma
_{j+1}^{\nu },\cdots ,\sigma _{l}^{\nu })\in \mathcal{C}_{l+1}(x)
\end{equation*}
\end{enumerate}

Noting $\mu =(\sigma _{1}^{\mu },\cdots ,\sigma _{l+1}^{\mu })$, we obtain
in both cases
\begin{equation}
\begin{array}{l}
(-1)^{l}\frac{e^{-ik\psi _{l}(x;\nu )}}{\left| \prod\limits_{j=1}^{l}\det
\left( \mathbb{I}+\mid \sigma _{j+1}^{\nu }-\sigma _{j}^{\nu }\mid \mathbb{P}%
_{j}^{\nu }\right) \right| ^{1/2}} \\
+(-1)^{l+1}\frac{e^{-ik\psi _{l+1}(x;\mu )}}{\left|
\prod\limits_{j=1}^{l+1}\det \left( \mathbb{I}+\mid \sigma _{j+1}^{\mu
}-\sigma _{j}^{\mu }\mid \mathbb{P}_{j}^{\mu }\right) \right| ^{1/2}}=0
\end{array}
\label{eq16}
\end{equation}
\end{lemma}

\textsc{Proof of Lemma }\ref{Lemma3}\textsc{\ :}\newline

In both case, using the strict convexity of compacts $(K_{j})_{j\in
\{1,\cdots ,N\}},$ we obtain $\mu \in \mathcal{C}_{l+1}(x)$.

\begin{remark}
In two cases, $\mu $ realize a transmission condition in $\sigma $, and thus
\begin{equation*}
\psi _{l}(x;\nu )=\psi _{l+1}(x;\mu )
\end{equation*}
\end{remark}

Now, we prove the formula (\ref{eq16}), in both case. In fact, we only have
to prove that
\begin{equation*}
{\left| \prod_{j=1}^{l}\det \left( \mathbb{I}+\mid \sigma _{j+1}^{\nu
}-\sigma _{j}^{\nu }\mid \mathbb{P}_{j}^{\nu }\right) \right| ^{1/2}}={%
\left| \left( \prod_{j=1}^{l+1}\det \left( \mathbb{I}+\mid \sigma
_{j+1}^{\mu }-\sigma _{j}^{\mu }\mid \mathbb{P}_{j}^{\mu }\right) \right)
\right| ^{1/2}}
\end{equation*}
Here $\sigma _{l+1}^{\nu }=\sigma _{l+2}^{\mu }=x$.\newline
In the first case, the proof is immediate.\newline
Under the hypothesis of the second case, we have
\begin{equation*}
\forall j\in \{1,\cdots ,j\}\ \ \mathbb{P}_{j}^{\nu }=\mathbb{P}_{j}^{\mu }
\end{equation*}
and, $\mu $ realize a transmission condition in $\sigma ,$ hence
\begin{equation*}
\mathbb{P}_{j+1}^{\mu }=S_{\mid \sigma _{j}^{\nu }-\sigma \mid }(\mathbb{P}%
_{j}^{\nu })=\mathbb{P}_{j}^{\nu }\left( \mathbb{I}+\mid \sigma _{j}^{\nu
}-\sigma \mid \mathbb{P}_{j}^{\nu }\right) ^{-1}.
\end{equation*}
Thus, we obtain
\begin{equation*}
\left( \mathbb{I}+\mid \sigma _{j+1}^{\nu }-\sigma \mid \mathbb{P}%
_{j+1}^{\mu }\right) \left( \mathbb{I}+\mid \sigma _{j}^{\nu }-\sigma \mid
\mathbb{P}_{j}^{\mu }\right) =\mathbb{I}+\mid \sigma _{j+1}^{\nu }-\sigma
_{j}^{\nu }\mid \mathbb{P}_{j}^{\nu }.
\end{equation*}
Taking determinant of previous formula, we get
\begin{equation*}
\det \left( \mathbb{I}+\mid \sigma _{j+1}^{\nu }-\sigma \mid \mathbb{P}%
_{j+1}^{\mu }\right) \det \left( \mathbb{I}+\mid \sigma _{j}^{\nu }-\sigma
\mid \mathbb{P}_{j}^{\mu }\right) =\det \left( \mathbb{I}+\mid \sigma
_{j+1}^{\nu }-\sigma _{j}^{\nu }\mid \mathbb{P}_{j}^{\nu }\right).
\end{equation*}
Moreover, we have
\begin{equation*}
S_{\mid \sigma _{j+2}^{\nu }-\sigma _{j+1}^{\mu }\mid }(\mathbb{P}%
_{j+1}^{\mu })=S_{\mid \sigma _{j+1}^{\nu }-\sigma \mid }\left( S_{\mid
\sigma _{j}^{\nu }-\sigma \mid }(\mathbb{P}_{j}^{\nu })\right)
\end{equation*}
As $\mu $ realize a transmission condition in $\sigma $ we obtain
\begin{equation*}
\mid \sigma _{j+1}^{\nu }-\sigma \mid +\mid \sigma _{j+1}^{\nu }-\sigma \mid
=\mid \sigma _{j+2}^{\nu }-\sigma _{j+1}^{\mu }\mid
\end{equation*}
thus
\begin{equation*}
S_{\mid \sigma _{j+2}^{\nu }-\sigma _{j+1}^{\mu }\mid }(\mathbb{P}%
_{j+1}^{\mu })=S_{\mid \sigma _{j+1}^{\nu }-\sigma _{j}^{\nu }\mid }(\mathbb{%
P}_{j}^{\nu }).
\end{equation*}
Then, we have
\begin{equation*}
\forall i\in \{j+1,\cdots ,l\}\ \ \mathbb{P}_{i}^{\nu }=\mathbb{P}%
_{i+1}^{\mu }
\end{equation*}
That's close proof of Lemma \ref{Lemma3}.$\Box $

\section{Proof of Theorem \ref{Theorem1}}

To proof this theorem, we apply stationary phase technics to the formula (%
\ref{eqD1}) and compare the result to the geometrical optic approximation.

\subsection{Stationary phase lemma}

Due to (\ref{eq7}) ($\mathbb{M}_{j}^{\nu }$ regular) we can apply the
iterative stationary phase lemma to $u_{l}^{D}(x)$ (see \cite{Cuvelier}) : $%
\forall x\in \Omega \setminus \mathcal{T}$
\begin{equation}
\begin{array}{c}
u_{l}^{D}(x)-u_{l-1}^{D}(x) \\
= \\
\left( \frac{ik}{4\pi }\right) ^{l}\left( \frac{2\pi }{k}\right)
^{l}\sum\limits_{\nu \in \mathcal{C}_{l}(x)}\frac{e^{i\frac{\pi }{4}%
\sum\limits_{j=1}^{l}\sgn \mathbb{M}_{j}^{\nu }}}{\left|
\prod\limits_{j=1}^{l}\det \mathbb{M}_{j}^{\nu }\right| ^{1/2}}\mathcal{D}%
_{l}(\nu )\frac{e^{-ik\psi _{l}(x;\nu )}}{\mid x-\sigma _{l}^{\nu }\mid }+O(%
\frac{1}{k}).
\end{array}
\label{eqD5}
\end{equation}
Thus, using definition of $\mathcal{D}_{l}(\nu )$ with $\nu \in \mathcal{C}%
_{l}(x),$ we obtain
\begin{equation*}
\mathcal{D}_{l}(\nu )=\left( 2^{l}\right) \left| \left\langle \xi ,n(\sigma
_{1}^{\nu })\right\rangle \right| \prod\limits_{j=1}^{l-1}{\frac{\left|
\left\langle {\frac{\sigma _{j+1}^{\nu }-\sigma _{j}^{\nu }}{\mid \sigma
_{j+1}^{\nu }-\sigma _{j}^{\nu }\mid },}n(\sigma _{j+1}^{\nu })\right\rangle
\right| }{\mid \sigma _{j+1}^{\nu }-\sigma _{j}^{\nu }\mid }}
\end{equation*}
Due to formula \ref{eq5} and \ref{eq6} (lemma \ref{Lemma2}), we have
\begin{equation*}
\sum\limits_{j=1}^{l}\sgn \mathbb{M}_{j}^{\nu }=2^{l}
\end{equation*}
and
\begin{eqnarray*}
\left| \prod\limits_{j=1}^{l}\det \mathbb{M}_{j}^{\nu }\right| ^{1/2}
&=&\prod\limits_{j=1}^{l}\frac{\left| \left\langle \xi _{j}^{\nu },n(\sigma
_{j}^{\nu })\right\rangle \right| }{\lambda _{j}^{\nu }}\left| \det (\mathbb{%
I}+\lambda _{j}^{\nu }\mathbb{P}_{j}^{\nu })\right| ^{1/2} \\
&=&\prod\limits_{j=1}^{l}\frac{\left| \left\langle \frac{\sigma _{j+1}^{\nu
}-\sigma _{j}^{\nu }}{\left| \sigma _{j+1}^{\nu }-\sigma _{j}^{\nu }\right| }%
,n(\sigma _{j}^{\nu })\right\rangle \right| }{\left| \sigma _{j+1}^{\nu
}-\sigma _{j}^{\nu }\right| }\left| \det (\mathbb{I}+\left| \sigma
_{j+1}^{\nu }-\sigma _{j}^{\nu }\right| \mathbb{P}_{j}^{\nu })\right| ^{1/2}
\end{eqnarray*}
with $\sigma _{l+1}^{\nu }=x.$ Taking into account that, for $\nu \in
\mathcal{C}_{l}(x),$ we have
\begin{eqnarray*}
\left| \left\langle \frac{\sigma _{2}^{\nu }-\sigma _{1}^{\nu }}{\left|
\sigma _{2}^{\nu }-\sigma _{1}^{\nu }\right| },n(\sigma _{1}^{\nu
})\right\rangle \right|  &=&\left| \left\langle \xi ,n(\sigma _{1}^{\nu
})\right\rangle \right|,\\
\left| \left\langle \frac{\sigma _{j+1}^{\nu }-\sigma _{j}^{\nu }}{\left|
\sigma _{j+1}^{\nu }-\sigma _{j}^{\nu }\right| },n(\sigma _{j}^{\nu
})\right\rangle \right|  &=&\left| \left\langle \frac{\sigma _{j}^{\nu
}-\sigma _{j-1}^{\nu }}{\left| \sigma _{j}^{\nu }-\sigma _{j-1}^{\nu
}\right| },n(\sigma _{j}^{\nu })\right\rangle \right| ,\,\;j=2,\ldots ,l,
\end{eqnarray*}
we find that
\begin{equation}
\begin{array}{c}
u_{l}^{D}(x)-u_{l-1}^{D}(x) \\
= \\
(-1)^{l}\sum\limits_{\nu =(\sigma _{1}^{\nu },\cdots ,\sigma _{l}^{\nu })\in
\mathcal{C}_{l}(x)}\frac{e^{-ik\psi _{l}(x;\nu )}}{\left| \left(
\prod\limits_{j=1}^{l-1}\det \left( \mathbb{I}+\mid \sigma _{j+1}^{\nu
}-\sigma _{j}^{\nu }\mid \mathbb{P}_{j}^{\nu }\right) \right) \det \left(
\mathbb{I}+\mid x-\sigma _{l}^{\nu }\mid \mathbb{P}_{l}^{\nu }\right)
\right| ^{1/2}} \\
+ \\
O(\frac{1}{k})
\end{array}
\label{eqD6}
\end{equation}

\subsection{Comparison with geometrical optic approximation}

To compare the previous formula with geometrical optic approximation given
by formula (\ref{GOA for Dirichlet}), we have to study the contributions of
the sets $\mathcal{C}_{l}(x)$ and $\mathcal{R}_{l}(x)$. We clearly have $%
\mathcal{R}_{l}(x)\subset \mathcal{C}_{l}(x)$. Using lemma \ref{Lemma3} we
obtain

\begin{remark}
The contributions of points in $\mathcal{C}(x)\setminus \mathcal{R}(x)$
cancel each other.
\end{remark}

Let $x\in \Omega \setminus \mathcal{T}$. Suppose that $\mathcal{C}%
_{l}(x)=\emptyset $ for $l>n$, we conclude the proof of Theorem \ref
{Theorem1} using the following remark

\begin{remark}
The only components of $\mathcal{C}(x)$ having a real contribution are :

\begin{itemize}
\item  all $\nu =(\sigma _{1}^{\nu })\in\mathcal{C}_1(x)$ \textit{%
\ coming through }$x$ such that $\nu$ realize a transmission
condition in $\sigma _{1}^{\nu }$ and
\begin{equation*}
\forall t>0\ \ \ \sigma _{1}^{\nu }-t\xi \in \Omega
\end{equation*}

\item  all $\nu =(\sigma _{1}^{\nu },\cdots ,\sigma
_{j}^{\nu })\in \mathcal{R}_j(x)$ \textit{\ coming through }$x$ ($j\leq n$).
\end{itemize}
\end{remark}
That's close proof of Theorem  \ref{Theorem1}.$\Box $
\section{Proof of Theorem \ref{Theorem2}}

To proof this theorem, we apply stationary phase technics to the formula (%
\ref{eqN1}) and compare the result to the geometrical optic approximation
(\ref{GOA for Neumann}).

Due to (\ref{eq7}) ($\mathbb{M}_{j}^{\nu }$ regular) we can apply the
iterative stationary phase lemma to $u_{l}^{N}(x)$(see \cite{Cuvelier}),we
obtain $\forall x\in \Omega \setminus \mathcal{T}$
\begin{equation}
\begin{array}{c}
u_{l}^{N}(x)-u_{l-1}^{N}(x) \\
= \\
\left( \frac{ik}{4\pi }\right) ^{l}\left( \frac{2\pi }{k}\right)
^{l}\sum\limits_{\nu \in \mathcal{C}_{l}(x)}\frac{e^{i\frac{\pi }{4}%
\sum\limits_{j=1}^{l}\sgn\mathbb{M}_{j}^{\nu }}}{\left|
\prod\limits_{j=1}^{l}\det \mathbb{M}_{j}^{\nu }\right| ^{1/2}}\mathcal{N}%
_{l}(\nu )\left\langle {\frac{x-\sigma _{l}^{\nu }}{\mid x-\sigma _{l}^{\nu
}\mid }},n(\sigma _{l}^{\nu })\right\rangle \frac{e^{-ik\psi _{l}(x;\nu )}}{%
\mid x-\sigma _{l}^{\nu }\mid }+O(\frac{1}{k})
\end{array}
\label{eqN5}
\end{equation}
Thus, we use (\ref{eq5}), (\ref{eq6}) and the definition of $\mathcal{N}%
_{l}$ to get :
\begin{equation}
\begin{array}{c}
u_{l}^{N}(x)-u_{l-1}^{N}(x) \\
= \\
\sum\limits_{\nu =(\sigma _{1}^{\nu },\cdots ,\sigma _{l}^{\nu })\in
\mathcal{C}_{l}(x)}\frac{e^{-ik\psi _{l}(x;\nu )}}{\left| \left(
\prod\limits_{j=1}^{l-1}\det \left( \mathbb{I}+\mid \sigma _{j+1}^{\nu
}-\sigma _{j}^{\nu }\mid \mathbb{P}_{j}^{\nu }\right) \right) \det \left(
\mathbb{I}+\mid x-\sigma _{l}^{\nu }\mid \mathbb{P}_{l}^{\nu }\right)
\right| ^{1/2}} \\
+ \\
O(\frac{1}{k})
\end{array}
\label{eqN6}
\end{equation}

To compare the previous formula with geometrical optic approximation given
by formula (\ref{GOA for Neumann}), we use previous results from the proof of theorem
\ref{Theorem1}.

\section{Conclusion}
We have proved the validity of the iterative Kirchhoff formulas (\ref{eqD1}) and (\ref{eqN1})
at high frequency. Numerical results for Dirichlet problem (\ref{Dirichlet}) can be founded in
\cite{Cuvelier}.

\bibliographystyle{alpha}
\bibliography{Dirichlet_Neumann}

\end{document}